\documentclass[a4paper,11pt]{article}
\usepackage{geometry,latexsym,amssymb,amsmath,color}
\usepackage[latin5]{inputenc}
\usepackage{enumerate}
\usepackage{enumitem}
\usepackage{hyperref}
\usepackage[all]{xy}
\usepackage{authblk}
\newtheorem{example}{Example}[section]
\newtheorem{Def}[example]{Definition}
\newtheorem{Exam}[example]{Example}

\newtheorem{Theo}[example]{Theorem}

\newtheorem{Cor}[example]{Corollary}
\newenvironment{Prf}{{\bf Proof:} } {\hfill $\Box$
\mbox{}}

\def\XMod{\operatorname{XMod}}
\def\Ker{\operatorname{Ker}}
\def\C{\operatorname{C}}
\def\DGpds{\operatorname{DGpds}}

\def\qed{\hfill $\Box$}
\def\Ker{\operatorname{Ker}}
\def\NSCM/(A,B,\alpha){\mathsf{NSCM/(A,B,\alpha)}}
\def\NSGGd/G{\mathsf{NSGGd/G}}
\def\DCatG{\mathsf{DCatG}}

\def\epsilon{\varepsilon}

\begin{document}
\title{\large\bf Normality and quotient in crossed modules over groupoids and double groupoids}
\author[1]{Osman MUCUK \thanks{Correspondence: mucuk@erciyes.edu.tr}}
\author[2]{Serap DEM\.{I}R \thanks{srpdmr89@gmail.com}}
\affil[1]{Department of Mathematics, Erciyes University Kayseri 38039, Turkey}
\affil[2]{Department of Mathematics, Erciyes University Kayseri 38039, Turkey}
\date{\vspace{-5ex}}
\maketitle

\noindent{\bf Key Words:} Quotient crossed module, double groupoid, quotient double groupoid
\\ {\bf Classification:} 20L05, 22A22, 18D35

\begin{abstract} We consider the categorical  equivalence  between  crossed modules over groupoids and
double groupoids with thin structures;  and by this equivalence, we prove how normality and quotient concepts are related in these two categories and give some examples of these objects.
\end{abstract}
\section{Introduction}
The concept of crossed module over groups introduced by Whitehead in \cite{23} and \cite{24} in the investigation of the properties of second relative homotopy groups for topological spaces, which  can be viewed as a 2-dimensional group \cite{2} has been widely used in: homotopy theory \cite{BHS}; the theory of identities among relations for group presentations \cite{4}; algebraic K-theory \cite{14}; and homological algebra \cite{13, 15}. See \cite [pp.49]{BHS}, for some discussion of the relation of crossed modules to crossed squares and so to homotopy 3-types.

On the one hand the categorical equivalence of  crossed modules over groups and group-groupoids,  which are internal groupoids in the category of groups and widely used in literature under the names \emph{2-groups}\cite{baez-lauda-2-groups}, {\em $\mathcal{G}$-groupoids} or \emph{group objects} in the category of groupoids  \cite{BS2},  was proved by Brown and Spencer in \cite[Theorem 1]{BS2}; and then some important results have been obtained by means of this equivalence. For example recently normal and quotient objects in these two categories have been compared and the corresponding objects in the category of group-groupoids have been characterized in \cite{19}.   The equivalence of these categories   has also been generalized by Porter in \cite[Section 3]{TP} to a more general  algebraic category $\C$, whose idea comes from Higgins \cite{Hig2} and Orzech \cite{Orz1,  Orz2} and called category of groups with operations. This result is used as a tool in the study  of  \cite{Ak-Al-Mu-Sa}. Applying Porter's result, the study of internal category theory in $\C$ was continued in the works of Datuashvili \cite{TD2} and \cite{TD3}.

On the other hand it was pointed out in \cite [Chapter 6]{BHS} that  the structure of crossed module is inadequate to give a proof of 2-dimensional Seifert-van-Kampen Theorem and hence one needs the idea of double groupoid which can be expressed as a groupoid object in the category of groupoids.  The categorical equivalence of crossed modules over groupoids and  double groupoids with thin structures was proved there.


In this paper as parallel  to the work studied in  \cite{19}  we define the normal subcrossed module and quotient crossed modules over groupoids; and then using the latter equivalence of the categories   we obtain the  normal subdouble groupoid and  quotient double groupoids. A motivation point of the relating  these objects  in both categories is to produce  more examples  of normal subdouble and quotient double  groupoids.
\section{Preliminaries on crossed modules over groupoids and double groupoids}
In the notations of  \cite [Chapter 6.1]{BHS},  a groupoid $G$ has a set $G_{0}$ of {\em objects},  a set $G_{1}$ of morphisms,  together with maps $\partial^{-}, \partial^{+}\colon G_{1}\rightarrow G_{0}$ and $\varepsilon  \colon G_{0} \rightarrow G_{1}$ such that $\partial^{-}\varepsilon=\partial^{+}\varepsilon=1_{G_{0}}$.
The maps $\partial^{-}$, $\partial^{+}$ are respectively called {\em initial} and {\em final} point maps  and the map $\varepsilon$ is called {\em object inclusion}.
If $a,b\in G$ and $\partial^{+}(a)=\partial^{-}(b)$, then the {\em composite} $a\circ b$ exists such that $\partial^{-}(a\circ b)=\partial^{-}(a)$ and $\partial^{+}(a\circ b)=\partial^{+}(b)$.
So there exists a partial composition defined by the map $ G\, {_{\partial^{+}}}\!\times_{\partial^{-}}G\rightarrow G, (a,b)\mapsto a\circ b$,
where $G\, _{\partial^{+}}\!\times_{\partial^{-}} G$ is the pullback of $\partial^{+}$ and $\partial^{-}$. Further, this partial composition is associative, for $x\in G_{0}$ the
element $\varepsilon (x)$ denoted by $1_x$ acts as the identity, and each morphism $a$ has an inverse $a^{-1}$
such that $\partial^{-}({a}^{-1})=\partial^{+}(a)$, $\partial^{+}({a}^{-1})=\partial^{-}(a)$, $a\circ {a}^{-1} =(\varepsilon \partial^{-})(a)$, $ {a}^{-1}\circ a=(\varepsilon \partial^{+})(a)$.
The map $G\rightarrow G$, $a\mapsto {a}^{-1}$,  is called the {\em inversion}.
    The set of all morphisms from $x$ to $x$ is a group, called \emph{object group}
at $x$, and denoted by $G(x)$. A groupoid $G$ is \emph{transitive (resp. totally intransitive)} if
$G(x,y)\neq\emptyset$ (resp. $G(x,y)=\emptyset$) for all $x,y\in G_0$ such that $x\neq y$.

A subgroupoid $H$ of $G$ is called {\em wide} if $H_0=G_0$ and a wide  subgroupoid $H$  is called {\em normal} if  $a\circ H(y)= H(x)\circ a$ for objects $x,y\in H_{0}$  and $a\in G(x,y)$. For example $\Ker f$, the kernel of a groupoid  morphism  $f\colon G\rightarrow K$,  is a normal subgroupoid of $G$.

Quotient groupoid is  formed as follows (see Higgins \cite[pp.86]{Hig1} and  \cite[pp.420]{8}). Let  $H$ be a normal subgroupoid of the groupoid  $G$.
The components of $H$ define a partition on $G_0$ and we write $[x]$ for the class containing $x$. Then $H$ also defines an equivalence relation on $G_1$ by $a\sim b$ for $a,b \in G$ if and only if $a=m\circ b \circ n$ for some $m,n\in H$. A partial composition $[a]\circ [b]$ on the morphisms is defined if and only if there exist $a_1\in [a]$, $b_1\in [b]$ such that $a_1\circ b_1$ is defined in $G$ and then $[a]\circ [b]=[a_1\circ b_1]$. This partial composition defines a groupoid on   classes $[x]$'s as objects. The groupoid defined in this manner is called {\em quotient groupoid}.

We remind the notion of double groupoids from  \cite [Chapter 6]{BHS}.

A {\em  double groupoid} denoted by $G=(G_2,G_1,G_0)$ has the sets $G_0$, $G_1$ and  $G_2$ of points or vertices, edges and squares respectively; and three structures of groupoids. The first one defined on $(G_1,G_0)$ has maps $\partial^{-}$,  $\partial^{+}$ and  $\varepsilon$ and  the composition denoted as multiplication. The other two are defined on $(G_2,G_1)$, a horizontal one  with maps $\partial_2^-$, $\partial_2^+$ and $\varepsilon_2$, composition denoted by $u+_2 v$; and a vertical one with maps $\partial_1^-$, $\partial_1^+$ and $\varepsilon_1$, composition denoted by  $u+_1 w$, satisfying some conditions given in details in   \cite [Chapter 6]{BHS}.

As a diagram a square $u\in G_2$ has bounding edges as
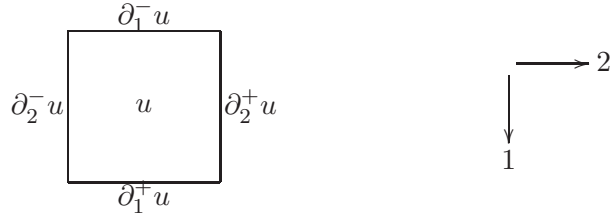
\[\begin{minipage}{0.5\textwidth}

\begin{xy}
*=<2cm,2cm>\txt{$u$}*\frm{-};
(14,0) *{\partial_2^+u} ;
(-14,0) *{\partial_2^-u} ;
(0,12) *{\partial_1^-u} ;
(0,-12) *{\partial_1^+u} ;
\end{xy}
\end{minipage}
\begin{minipage}{0.5\textwidth}

\begin{displaymath}
\xymatrix{
\ar[d]  \ar[r] & 2 \\
1 & }
\end{displaymath}
\end{minipage}\]
For an edge $a\in G_1$,  the identities $\varepsilon_1(a)$ and $\varepsilon_2(a)$ for vertical and horizontal composites of squares which are called {\em degeneracies} have respectively the boundaries

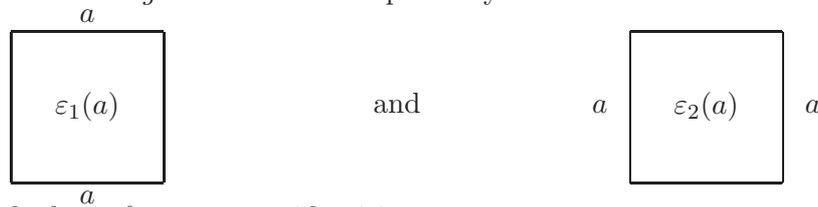
\begin{minipage}{0.5\textwidth}
\[\begin{xy}
*=<2cm,2cm>\txt{$\epsilon_1(a)$}*\frm{-};
(14,0) *{};
(-14,0) *{} ;
(0,12) *{a} ;
(0,-12) *{a} ;
\end{xy}\]
\end{minipage}
and
\begin{minipage}{0.5\textwidth}
\[\begin{xy}
*=<2cm,2cm>\txt{$\epsilon_2(a)$}*\frm{-};
(14,0) *{a} ;
(-14,0) *{a} ;
(0,12) *{} ;
(0,-12) *{} ;
\end{xy}\]
\end{minipage}
where the unlabeled edges of squares are identities.

A \emph{morphism} $f\colon G\rightarrow K$ of double groupoids is a triple of  functions
\[f_{i}\colon G_{i}\rightarrow K_{i},~~~~~(i=0,1,2) \]
which commute with all structural maps of three groupoids (faces, degeneracies, compositions etc.).  Then  we have a category $\DCatG$ of double groupoids and morphisms between them.

The following example of double goupoid  is from  \cite[Definition 6.1.10]{BHS}.
\begin{Exam} \label{2shellDgpd}{\em For a groupoid $G$, a  double groupoid of `squares' or 2-`shells' in $G$ denoted by $\square'G$ is defined as follows:

The groupoid structure of  $(\square' G)_1$ is the same as that of $G$. The squares of  $(\square' G)_2$ consist of quadruples $\left ( a ~\begin{array}{c} c \\   b \\ \end{array} ~ d\right)$ for $a,b,c,d\in G_1$;   the horizontal and vertical compositions are defined by
\[\left ( a ~\begin{array}{c} c \\   b \\ \end{array} ~ d\right)+_2\left(d ~\begin{array}{c} g \\  f \\ \end{array} ~ h\right)=\left( a ~\begin{array}{c} cg \\   bf\\ \end{array} ~ h\right)\]

\[\left ( a ~\begin{array}{c} c \\   b \\ \end{array} ~ d\right)+_1 \left( f ~\begin{array}{c} b \\   g \\ \end{array} ~ h\right)=\left( af ~\begin{array}{c} c \\   g \\ \end{array} ~ dh\right). \]

The double groupoid of  commutative squares or commutative 2-shells is denoted by $\square G$. The squares of  $(\square G)_2$ consist of quadruples $\left ( a ~\begin{array}{c} c \\   b \\ \end{array} ~ d\right)$ for $a,b,c,d\in G_1$ such that $ab=cd$.
}\end{Exam}

\begin{Def} {\em \cite[Definition 6.4.1]{BHS}
A \emph{thin structure} for a double groupoid $G$ is a morphism  of double groupoids \[\circleddash \colon \square (G_1)\rightarrow G \] which is identity on $G_0$ and $G_1$.  A  2-dimensional element $\circleddash(\alpha)\in G_2$ for an $\alpha\in (\square G_1)_2 $ is called a {\em thin square}.}\end{Def}

  In  \cite [Section 6.4]{BHS} a criterion for  the existence of a thin  structure is  given in the terms of  thin squares.

Let $\DGpds$ be the category of double groupoids with thin structures and morphisms between them preserving the thin structures.

\begin{Exam}  {\em  For a groupoid $G$ one can easily check that , the double groupoid $\square' G$ of 2-shells defined in Definition \ref{2shellDgpd} has a trivial double groupoid morphism  as a thin structure  \[\circleddash \colon \square (G_1)\rightarrow \square' G \]
mapping the commutative 2-shells to themselves.} \end{Exam}

     We know from  \cite[Theorem 6.4.11]{BHS} that in a double groupoid which has a thin structure, the horizontal and vertical groupoid structures  on squares  are isomorphic. Hence  the conditions of  being subgroupoid and normal subgroupoid in one direction  can be transferred to the other one and therefore we state and  prove the theorems in Section \ref{Sectionnormaldgpd} for only horizontal composite $+_2$.

In  \cite[Definition 6.2.1]{BHS}  a \emph{crossed module over a groupoid} $P=(P_1,P_0)$ is defined as a morphism of groupoids \[\mu\colon M\rightarrow P\]
in which $M=(M_1,P_0)$ is a totally intransitive  groupoid and  $\mu$ is identity on objects. Here  $+$ is used  for the composition in groupoid  $M$.  The groupoid $P$ operates on the right on $M$ and the action  denoted by \[M(p)\times P(p,q)\rightarrow M(q)\]\[(x,a)\mapsto x^{a}\]
satisfies the usual axioms of an action (i) $x^{1}=x$, $(x^{a})^{b}= x^{ab}$ and (ii) $(x+y)^a=x^{a}+y^{a}$ and the following conditions are satisfied:
\begin{enumerate}[label=\textbf{\roman{*}.}, leftmargin=1cm]
\item[CM1.] $\mu$ preserves the action, i.e.,  $\mu(x^a)=\mu (x)^a=a^{-1}\mu(x)a$.

\item[CM2.] For all $c\in M{(p)}$, $\mu(c)$ acts on $M$ by conjugation \[x^{\mu(c)}=-c+x+c.  \]
for any $x\in M(p)$.
\end{enumerate}

\begin{Exam} {\em If $H$ is a totally intransitive normal subgroupoid of $G$,  then the inclusion $\imath\colon H\rightarrow G$ is a crossed module  over groupoids with the action of conjugation. }\end{Exam}

For a  topological example called {\em fundamental crossed module} over groupoids see \cite[Example 6.2.2]{BHS}.

A \emph{morphism} $(f_2,f_1)\colon (\mu \colon M\rightarrow P) \rightarrow ( \nu\colon N\rightarrow Q)$ of crossed modules over groupoids   is a pair of groupoid morphisms $f_2\colon M\rightarrow N$ and $f_1\colon P \rightarrow Q$ with the same induced map of vertices  such that $f_1\mu=\nu f_2$ and $f_2(m^{a})=f_2(m)^{f_1(a)}$ for $m\in M$ and $a\in P$. Hence crossed modules over groupoids and morphisms of them constitute a category, denoted by $\XMod$.  From now on by a crossed module we mean a crossed module over groupoids.

The following result was proved in \cite[Section 6.6]{BHS}. Since we  need some details of  the proof in the following section, we give a sketch proof.

\begin{Theo}\label{teo1}
The category $\XMod$ and $\DGpds$ are equivalent.
\end{Theo}
\begin{Prf}
A functor $\lambda\colon \XMod\rightarrow \DGpds$ is defined as follows:  For a crossed module $\mu \colon M\rightarrow P$, there is an  associated double groupoid $G=\lambda(\mu \colon M\rightarrow P)$ in which   $G_{0}=P_0$, $G_{1}=P_1$ and  the set $G_2$ of squares consists of the elements
\[ \left(m ~;~\left( a ~\begin{array}{c} c \\   b \\ \end{array} ~ d\right)\right)\]
such that $ a, b, c, d\in P_1$ and $\mu(m)=b^{-1}a^{-1}cd$.
The elements of $G_2$ can be represented by
\[\begin{xy}
*=<2cm,2cm>\txt{$m$}*\frm{-};
(14,0) *{d} ;
(-14,0) *{a} ;
(0,12) *{c} ;
(0,-12) *{b} ;
\end{xy}\]

The horizontal and vertical compositions of squares are respectively defined by
\[\left(m ~;~\left( a ~\begin{array}{c} c \\   b \\ \end{array} ~ d\right)\right)+_2\left(n ~;~\left( d ~\begin{array}{c} g \\   f \\ \end{array} ~ h\right)\right)=\left(m^f+n ~;~\left( a ~\begin{array}{c} cg \\   bf \\ \end{array} ~ h\right)\right)\]

\[ \left(m ~;~\left( a ~\begin{array}{c} c \\   b \\ \end{array} ~ d\right)\right) +_1\left(n ~;~\left( f ~\begin{array}{c} b \\   g\\ \end{array} ~ h\right)\right)=\left(n+m^h ~;~\left( af ~\begin{array}{c} c \\  g \\ \end{array} ~ dh\right)\right)\]

A thin  structure
$\circleddash \colon \square G_1 \rightarrow G$ on $G$ is defined by
\[\left ( a ~\begin{array}{c} c \\   b \\ \end{array} ~ d\right)\mapsto\left (1 ~;~ \left(a ~\begin{array}{c} c \\  b \\ \end{array} ~ d\right)\right).\]

Conversely for a double groupoid  $G$, an associated crossed module  $\gamma (G)=(\mu\colon M \rightarrow P )$  is given by  $P=(G_1,G_{0})$ and  $M=(M_1,G_0)$ in which $M_1$ has the squares $u\in G_2$ in the form of   \[ \begin{xy}
*=<2cm,2cm>\txt{$u$}*\frm{-};
(14,0) *{} ;
(-14,0) *{} ;
(0,12) *{\partial_1^{-}u} ;
(0,-12) *{} ;
\end{xy}.\]

Here the boundary  morphism $\mu$ is  ${\partial_1}^{-}\colon M_1\rightarrow P_1$ and the action $u^a $ is defined by
\begin{align*}
&\begin{xy}
*=<2cm,2cm>\txt{$u^a$}*\frm{-};
(14,0) *{} ;
(-14,0) *{} ;
(0,12) *{a^{-1}\partial_1^{-}(u) a} ;
(0,-12) *{} ;
\end{xy}
\end{align*}
for $a\in G_1(x,y)$ and $u\in M_1(x)$;  and the conditions CM1 and CM2 are satisfied.  Hence we have a functor $\gamma\colon \DGpds\rightarrow \XMod$.

The functors $\gamma\lambda\colon \XMod\rightarrow \XMod$ and $\lambda\gamma\colon \DGpds\rightarrow \DGpds$ are   naturally equivalent to the identity functors $1_{\XMod}$ and $1_{\DGpds}$ respectively.
\end{Prf}

We remark that as related to Theorem \ref{teo1},  the  categorical  equivalence  of  the crossed modules over groups and special double groupoids with connections  was proved in  \cite{BS1} and using this equivalence, in  \cite{SD}  the normal and quotient  objects in these categories  are characterized. In the next section this  is generalized to the equivalence given in Theorem \ref{teo1}.

\section{Normality and quotient in double groupoids} \label{Sectionnormaldgpd}

In this section using the equivalence of the categories as stated in Theorem \ref{teo1},  we will determine the normal and quotient objects in the category  $\DGpds$ of double groupoids with thin structures. Before this determination we initially  define  subcrossed module and normal subcrossed module over groupoids as similar to the group cases given in \cite{20}.
\begin{Def}  \label{scm}
\em{  A crossed module $\nu \colon N\rightarrow Q$ is called a \emph{subcrossed module} of the crossed module $\mu \colon M\rightarrow P$ if the followings hold:
\begin{enumerate}
\item[SCM1.] $N$ is a subgroupoid of $M$ and $Q$ is a subgroupoid of $P$,
\item[SCM2.] $\nu$ is restriction of $\mu$ to $N$ and
\item[SCM3.] The action of $Q$ on $N$ is the restriction of the action of $P$ on $M$.
\end{enumerate}}
\qed
\end{Def}
\begin{Def}\label{ncm}
\em{  A subcrossed module $\nu \colon N\rightarrow Q$ of the crossed module $\mu \colon M\rightarrow P$ is called \emph{normal } if
\begin{enumerate}
\item [NCM1.] $Q$ is a normal subgroupoid of $P$,
\item [NCM2.] $n^{a}\in N(q)$ for all $a\in P(p,q)$, $n\in N(p)$ and
\item [NCM3.] $-m+m^{b}\in N(p)$ for all $b\in Q(p)$ and $m\in M(p)$.
\end{enumerate}}
\qed
\end{Def}

Remark that in this definition $N$ also becomes a normal subgroupoid of $M$. Indeed for $n \in N(p)$ and $m\in M(p)$, we have  $n^{\mu(m)}\in N(p)$ by condition [NCM2] and  $n^{\mu(m)}=-m+n+m  \in N(p)$ by [CM2].

\begin{Exam}\label{ker}
\em{The {\em  kernel} of a crossed module morphism $ (f_2,f_1)\colon (\mu \colon M\rightarrow P)\rightarrow (\nu \colon N\rightarrow Q) $ defined by
\[\mu^{\star}\colon \Ker f_2\rightarrow \Ker f_1\]
as  a restriction of $\mu$ is a normal subcrossed module of the crossed module $\mu \colon M\rightarrow P$.}
\end{Exam}

\begin{Theo} \label{Quotientcrospond}
Let $\nu \colon N\rightarrow Q$ be a normal subcrossed module of $\mu \colon M\rightarrow P$. If the groupoid $Q$ is totally intransitive, then $\rho \colon M/N\rightarrow P/Q$ is a crossed module over groupoids.
\end{Theo}
\begin{Prf}
By the construction of quotient groupoid we know that when  $Q$ is totally intransitive, object set of the quotient groupoid $P/Q$ is $P_0$.   Similarly since $N$ is totally intransitive, $(M/N)_0=M_0=P_0$. Hence  $M/N$ and $P/Q$ have the same objects. Since $M$ and $N$ are totally intransitive, so also $M/N$ is. Therefore the induced  morphism $\rho\colon M/N\rightarrow P/Q $ is well defined.

An action  of $P/Q$ on $M/N$ is defined  by
\[ (M/N)(p)\times (P/Q)(p,q)\rightarrow (M/N)(q)\]
\[ (N(p)+m,~Q(p)a)\mapsto (N(p)+m)^{Q(p)a} =N(q)+m^a\]
for $a\in P(p,q)$ and $m\in M(p)$. We now prove that this  action is  well defined.  If  $m_1\in (N(p)+m)$ and $a_1\in Q(p)a$, there exist $n\in N(p) $ and $b\in Q(p)$ such that $m_1=n+m$ and $a_1=ba$. Here  $-m+m^b=n_1 \in N(p)$ by [NCM3] and hence $m^b=m+n_1$. There exists an $n_2\in N(p)$ such that $m+n_1=n_2+m$ by [NCM2].  By the action of $Q$ on $N$ and [NCM2] we have $n^{ba}+n_2^{a}\in N(q)$. If we write  $n_3$ for $n^{ba}+n_2^{a}$, then we get
\begin{align*}
m_1^{a_1} &=(n+m)^{ba} =n^{ba}+m^{ba} =n^{ba}+(m^b)^a\\
&=n^{ba}+(m+n_1)^{a}=n^{ba}+(n_2+m)^{a}\\
&=n^{ba}+n_2^{a}+m^{a}\\
&=n_3+m^a.
\end{align*}
Hence  $n_3+m^a\in (N(q)+m^a)$, i.e., the  action is well defined. It can be easily checked  that  the conditions [CM1] and [CM2] are satisfied and hence $\rho\colon M/N\rightarrow P/Q $  becomes a crossed module over groupoids.
\end{Prf}

We call  $\rho \colon M/N\rightarrow P/Q$  as  {\em quotient crossed module} over groupoids.

Let $f\colon P\rightarrow Q$ be a morphism of groupoids.  To use in the following example we remind  from \cite[Proposition 26]{Hig1} that   if  $f$ is injective on objects, then $\Ker f$ is totally intransitive.
\begin{Exam}\label{quotintdgpdexam}
\em{Let $\mu^{\star}\colon \Ker f_2\rightarrow \Ker f_1$ be  the kernel of a crossed module morphism $ (f_2,f_1)\colon (\mu \colon M\rightarrow P)\rightarrow (\nu \colon N\rightarrow Q) $  as defined in Example \ref{ker} such that $f_1$ is injective on objects. Then by Theorem \ref{Quotientcrospond} the induced  morphism $\rho\colon M/\Ker f_2\rightarrow P/\Ker f_1$ is a quotient crossed module.}
\end{Exam}

Since in a double groupoid with a thin structure the horizontal and vertical groupoid structures  on squares  are isomorphic, in the following theorems we use only horizontal composite $+_2$.

\begin{Theo} \label{direct}
Let $\nu \colon N\rightarrow Q$ be a subcrossed module of a crossed module $\mu \colon M\rightarrow P$. Suppose that  $H=(H_{2},H_{1},H_{0})$ and $G=(G_{2},G_{1},G_{0})$ are respectively
the double groupoids corresponding to these crossed modules. Then the followings are satisfied.
\begin{enumerate}
\item $(H_1,H_{0})$ is a subgroupoid of $(G_1,G_{0})$.

\item $(H_2,H_{1},+_2)$ is a subgroupoid of $(G_2,G_{1},+_2)$.

\item The thin structure $\circleddash\colon \square (H_1)\rightarrow H$ on $H$ is a restriction of that $\circleddash\colon \square (G_1)\rightarrow G$ on $G$.
\end{enumerate}
\end{Theo}

\begin{Prf} \begin{enumerate}
\item By the proof of  Theorem \ref{teo1},  we know that $(H_{1},H_0)=Q$ and $(G_{1},G_0)=P$; and   $(H_1,H_{0})$ becomes a  subgroupoid of $(G_1,G_{0})$ by [SCM1] of Definition \ref{scm}.

\item  $H_2$ consists of the elements
\[\alpha=\left(m ~;~ \left(a ~\begin{array}{c} c \\   b \\ \end{array} ~ d\right)\right)\]
for  $m\in N$,  $a,b,c,d\in Q_{1}$ such that  $\nu(m)=b^{-1}a^{-1}cd$; and   such a quintuple has an inverse
\[-_2\alpha=\left((-m)^{b^{-1}} ~;~ d ~\begin{array}{c} c^{-1} \\   b^{-1} \\ \end{array} ~ a\right)\]
for horizontial composite. Since $N$ is a groupoid and $P$ acts on $N$, we have that  $(-m)^{b^{-1}}\in N$. Hence $-_2\alpha\in H_2$ and thus   $(H_{2},H_{1},+_2)$ is a subgroupoid of $(G_{2},G_{1},+_2)$.

\item  This is obvious by the  details of the proof of Theorem  \ref{teo1}.
\end{enumerate}
\end{Prf}

Hence we can state the definition of a subdouble groupoid as follows.
\begin{Def}
{\em A double groupoid $H=(H_{2},H_{1},H_{0})$ is called \emph{a subdouble groupoid} of $G=(G_{2},G_{1},G_{0})$ if the followings are satisfied:
\begin{enumerate}
\item[SDG1.] $(H_1,H_{0})$ is a subgroupoid of $(G_1,G_{0})$.

\item[SDG2.] $(H_2,H_{1},+_2)$ is a subgroupoid of $(G_2,G_{1},+_2)$.

\item[SDG3.] The thin structure $\circleddash\colon \square (H_1)\rightarrow H$ on $H$ is a restriction of the thin structure $\circleddash\colon \square (G_1)\rightarrow G$ on $G$.
\end{enumerate}}
\end{Def}

We now prove that the converse of Theorem \ref{direct} is also satisfied.
\begin{Theo}\label{convers} Let $H=(H_{2},H_{1},H_{0})$ be a subdouble groupoid of a double groupoid $G=(G_{2},G_{1},G_{0})$.  Assume that   $\nu\colon N\rightarrow Q$ and $\mu\colon M\rightarrow P$ are respectively  the crossed modules corresponding to $H$ and $G$. Then the crossed module
$\nu\colon N\rightarrow Q$ is a subcrossed module of the crossed module $\mu\colon M\rightarrow P$.
\end{Theo}
\begin{Prf}We need to prove that the conditions of Definition \ref{scm} are satisfied.

SCM1. By  Theorem \ref{teo1},  $Q=(H_{1},H_{0})$, $P=(G_{1},G_{0})$ and by [SDG1] $(H_1,H_{0})$ is subgroupoid of $(G_1,G_{0})$. Hence  $Q$ is subgroupoid of $P$. Further  $N$ and $M$  consist of respectively specific squares of $H_{2}$ and $G_{2}$  as stated in the proof of Theorem \ref{teo1} and by [SDG2] $(H_{2},H_{1},+_2)$ is a subgroupoid of  $(G_{2},G_{1},+_2)$. Hence $N$ becomes a subgroupoid of $M$.

SCM2. By the details of the proof of Theorem  \ref{teo1}, the crossed modules $\nu\colon N\rightarrow Q$ and $\mu\colon M\rightarrow P$ are defined by

\begin{center}
$\nu\colon N\rightarrow Q, n\mapsto \partial_1^{-}n$\\
$\mu\colon M\rightarrow P, m\mapsto \partial_1^{-}m$
\end{center}
and hence  $\nu$ is a restriction of $\mu$.

SCM3. It is clear that the action of $P$ on  $N$ is a restriction of the action of $P$ on $M$.

Hence the crossed module $\nu\colon N\rightarrow Q$ becomes a subcrossed module of the crossed module $\mu\colon M\rightarrow P$ as required.
\end{Prf}

As a result of Theorem \ref{direct} and Theorem \ref{convers} the following corollary can be stated.
\begin{Cor}
Let $G$ be a double groupoid and $\mu\colon M\rightarrow P$ the  crossed module corresponding to $G$.
Then the category of the  subdouble groupoids of $G$ and the category of  subcrossed modules of  $\mu\colon M\rightarrow P$ are equivalent.
\end{Cor}

As we can see from the proof of Theorem \ref{Theodirec} in the subdouble groupoid corresponding to a normal subcrossed module the horizontal and the vertical subgroupoids are not necessarily wide. So we need the definition of a normal subgroupoid without wideness condition  as follows.

\begin{Def}\label{ansg}
\emph{Let $G$ be a groupoid and $H$ be a subgroupoid of $G$. Then $H$ is called a \emph{non-wide normal subgroupoid} of $G$ if $gag^{-1}\in H(x)$
when $g\in G(x,y)$ and $a\in H(y)$ for $x,y\in H_0$.}
\end{Def}

\begin{Theo}\label{Theodirec}
Let $\nu\colon N\rightarrow Q$ be a normal subcrossed module of a crossed module  $\mu\colon M\rightarrow P$. Suppose that $H=(H_{2},H_{1},H_{0})$ and $G=(G_{2},G_{1},G_{0})$ are respectively the double groupoids corresponding  to these crossed modules. Then we have the followings.
\begin{enumerate}
\item $(H_{1},H_{0})$ is a normal subgroupoid of $(G_{1},G_{0})$.

\item The groupoid $(H_{2},H_{1},+_2)$ is a non-wide normal subgroupoid of $(G_{2},G_{1},+_2)$.

\item The thin structure $\circleddash\colon \square (H_1)\rightarrow H$ on $H$ is the restriction of the thin structure $\circleddash\colon \square (G_1)\rightarrow G$ on $G$.
\end{enumerate}
\end{Theo}
\begin{Prf}
\begin{enumerate}
\item  By the equivalence given in Theorem \ref{teo1}, $(H_{1},H_{0})=Q$, $(G_{1},G_{0})=P$ and by [NCM1] in Definition \ref{ncm}  $Q$ is normal subgroupoid of $P$. Hence  $(H_{1},H_{0})$ is a normal subgroupoid of $(G_{1},G_{0})$.

\item   Let
$\alpha=\left(m ~;~\left( a ~\begin{array}{c} c \\   b \\ \end{array} ~ d\right)\right)\in G_{2}$ and
$\beta=\left(n ~;~ \left(d ~\begin{array}{c} g \\   f \\ \end{array} ~ d\right)\right)\in H_{2}$. Then we have that
\begin{align*}
\alpha+_2\beta-_2\alpha&=\left(m ~;~\left( a ~\begin{array}{c} c \\   b \\ \end{array} ~ d\right)\right)+_2\left(n ~;~ \left(d ~\begin{array}{c} g \\   f \\ \end{array} ~ d\right)\right)+_2\left((-m)^{b^{-1}} ~;~\left( d ~\begin{array}{c} c^{-1} \\   b^{-1} \\ \end{array} ~ a\right)\right)\\
&=\left(m^f+n ~;~\left( a ~\begin{array}{c}  cg\\  bf  \\ \end{array} ~ d\right)\right)+_2\left((-m)^{b^{-1}} ~;~\left( d ~\begin{array}{c}  c^{-1} \\  b^{-1}  \\ \end{array} ~ a\right)\right) \\
&=\left((m^f+n)^{b^{-1}}+(-m)^{b^{-1}} ~;~\left( a ~\begin{array}{c} cgc^{-1} \\   bfb^{-1} \\ \end{array} ~ a\right)\right).
\end{align*}
By [NCM1] since $Q$ is normal subgroupoid of $P$ we have that  $cgc^{-1}\in Q$ and $bfb^{-1}\in Q$ for $f,g\in Q$ and $c,b\in P$. Here  we need to prove that $(m^f+ n)^{b^{-1}}+(-m)^{b^{-1}}\in N$. If we write $n_1$ for  $-m+m^f\in N$ by [NCM3], then we have that
\begin{align*}
(m^f+ n)^{b^{-1}}+(-m)^{b^{-1}}&=(m^f+ n-m)^{b^{-1}}\\
&=(m-m+m^f-m+m+ n-m)^{b^{-1}}\\
&=(m+n_1-m+m+ n-m)^{b^{-1}}
\end{align*}

Here by [NCM2]  $m+n_1-m$, $m+ n-m\in N$ and for $b\in Q$, $n,n_1\in N$ and $m\in M$; $(m^f+ n)^{b^{-1}}+(-m)^{b^{-1}}\in N$. We now  show that the boundary condition
\[\mu((m^f+ n)^{b^{-1}}+(-m)^{b^{-1}})=bf^{-1}b^{-1}a^{-1}cgc^{-1}a\] is satisfied.
Since $\mu(m)=b^{-1}a^{-1}cd$, $\mu(n)=f^{-1}d^{-1}gd$ we have that
\begin{align*}
\mu((m^f+ n)^{b^{-1}}+(-m)^{b^{-1}})&=\mu(m^{(fb^{-1})})\mu(n^{b^{-1}})\mu((-m)^{b^{-1}})\\
&=(bf^{-1}\mu(m)fb^{-1})(b\mu(n)b^{-1})(b\mu(-m)b^{-1}) \tag{by [CM1]}\\
&=bf^{-1}b^{-1}a^{-1}cgc^{-1}a.
\end{align*}

Hence we obtain that $\alpha+_2\beta-_2\alpha\in H_{2}$, i.e., $(H_2,H_1,+_2)$
is non-wide normal subgroupoid of $(G_2,G_1,+_2)$ in the sense of Definition \ref{ansg}.

\item  This comes from the fact that   $H$ is a subdouble groupoid of $G$.
\end{enumerate}
\end{Prf}

Hence we can give the definition of a normal subdouble groupoid as follows.

\begin{Def}\label{ndg}
{\em A subdouble groupoid $H=(H_{2},H_{1},H_{0})$  of a double groupoid $G=(G_{2},G_{1},G_{0})$ is  called \emph{normal subdouble groupoid} if the followings hold:

\begin{enumerate}
\item[NDG1.] $(H_{1},H_{0})$ is a normal subgroupoid of $(G_{1},G_{0})$.

\item[NDG2.]The groupoid $(H_{2},H_{1},+_2)$ is a non-wide normal subgroupoid of $(G_{2},G_{1},+_2)$.
\end{enumerate}}
\end{Def}

We use the equivalence of the categories given in Theorem \ref{Theodirec} to produce the following examples of normal subdouble groupoids.

\begin{Exam}{\em  Let $(f_{2},f_{1})\colon (\mu\colon M\rightarrow P)\rightarrow (\nu\colon N\rightarrow Q)$ be a  morphism of crossed module with the kernel  $\mu^{\star}\colon Ker f_{2}\rightarrow Kerf_{1}$ as defined in Example \ref{ker}, which is a normal subcrossed module of $\mu\colon M\rightarrow P$.  Let $H=(H_2,H_1,H_0)$ be the double groupoid corresponding to $\mu^{\star}\colon Ker f_{2}\rightarrow Kerf_{1}$, where $H_0=P_0$, $H_1=\Ker f_1$ and $H_2$ consists of  \[ \left(n ~;~\left( a ~\begin{array}{c} c \\   b \\ \end{array} ~ d\right)\right)\]
such that $ a, b, c, d\in \Ker f_1$ and $n\in \Ker f_2$ with $\mu^{\star}=b^{-1}a^{-1}cd$.  Assume that $G=(G_2,G_1,G_0)$ is the double groupoid corresponding to $\mu$ and hence defined by $G_0=P_0$, $G_1=P_1$ and $G_2$ consists of  \[ \left(m ~;~\left( a ~\begin{array}{c} c \\   b \\ \end{array} ~ d\right)\right)\]
such that $ a, b, c, d\in P_1$ and $m\in M_1$ with $\mu(m)=b^{-1}a^{-1}cd$. Then by Theorem  \ref{Theodirec}  $H$ becomes a normal subdouble groupoid of $G$.} \end{Exam}

\begin{Exam} {\em  Let $H$ be a normal subgroupoid of  $G$. Then 2-shell groupoid $\square' H$ becomes a normal subdouble groupoid of the double groupoid  $\square' G$}.\end{Exam}

\begin{Theo} \label{TheoConver}
Let $H=(H_{2},H_{1},H_{0})$ be a  normal subdouble groupoid of a double groupoid $G=(G_{2},G_{1},G_{0})$.  Let $\nu\colon N\rightarrow Q$ and $\mu\colon M\rightarrow P$ be respectively  the crossed modules corresponding  to $H$ and $G$. Then the crossed module
$\nu\colon N\rightarrow Q$ is a normal subcrossed module of the crossed module $\mu\colon M\rightarrow P$.
\end{Theo}
\begin{Prf} We must show that the conditions of Definition \ref{ncm} are satisfied:

NCM1. Since $Q=(H_{1},H_0)$, $P=(G_{1},G_0)$ and by [NDG1] $(H_{1},H_0)$ is a normal subgroupoid of $(G_{1},G_0)$; $Q$ becomes a normal subgroupoid of $P$.

NCM2.  We need to prove that for $a\in P(p,q)$ and  $n\in N(p)$, $n^{a}$ is in $N(q)$. Here $ n^{a}$ is of the form
\[\begin{xy}
*=<2cm,2cm>\txt{$n^{a}$}*\frm{-};
(14,0) *{} ;
(-14,0) *{} ;
(0,12) *{a^{-1}\nu(n)a} ;
(0,-12) *{} ;
\end{xy}.\]

By [NDG1] since $(H_{1},H_0)$ is a normal subgroupoid of $(G_{1},G_0)$, it implies  that $a^{-1}\nu(n)a\in H_{1}$ and hence $n^{a}\in N(q)$.

NCM3.  We prove that  $-m+m^{b}\in N$ for  $b\in Q=H_{1}$ and $m\in M$. Here the square
$-m+m^{b}$ can be denoted by a square as follows;

\[\begin{minipage}{0.25\textwidth}
\[\begin{xy}
*=<2cm,2cm>\txt{$-m+m^{b}$}*\frm{-};
(14,0) *{} ;
(-14,0) *{} ;
(0,12) *{\mu(-m)b^{-1}\mu(m)b} ;
(0,-12) *{} ;
\end{xy}\]
\end{minipage}\]

Here $b\in Q=H_{1}$,~ $\mu(m)\in P=G_{1}$ and by [NDG1] since  $(H_{1},H_0)$ is normal subgroupoid of $(G_{1},G_0)$, we have that $\mu(-m)b^{-1}\mu(m)\in Q$ and $\mu(-m)b^{-1}\mu(m)b\in Q$. As a result we get  \[-m+m^{b}\in N\]

Thus $\nu\colon N\rightarrow P$ is a normal subcrossed module of $\mu\colon M\rightarrow P$.
\end{Prf}

Hence the following corollary is stated as a result of Theorem \ref{Theodirec} and Theorem \ref{TheoConver}.

\begin{Cor}\label{EquivalentofNormals}
Let $G$ be a double groupoid and $\mu\colon M\rightarrow P$ the crossed module corresponding to $G$. Then  the category of the normal subdouble groupoids of $G$ and the category of normal subcrossed modules of the crossed module $\mu\colon M\rightarrow P$ are equivalent.
\end{Cor}

\begin{Def} {\em  Let   $\nu\colon N\rightarrow Q$ be a  normal subcrossed module of $\mu\colon M\rightarrow P$ such that $Q$ is totally intransitive. Then the double groupoid corresponding to the  quotient crossed module $\rho\colon M/N\rightarrow P/Q$ by Theorem \ref{Quotientcrospond}  is called {\em quotient double groupoid}}.\end{Def}

Hence we can characterize a quotient groupoid as follows: Let  $H=(H_2,H_1,H_0)$ be a  normal subdouble groupoid of $G=(G_2,G_1,G_0)$ such that the groupoid  $(H_1,H_0)$ is totally intransitive. Then by Corollary \ref{EquivalentofNormals} we can construct the quotient double groupoid $G/H$  whose points are same as that of  $G$, edges are the morphisms of  quotient groupoid $G_1/H_1$, and the squares are in the form of

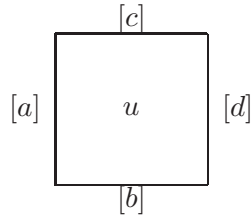
\[\begin{minipage}{0.5\textwidth}
\begin{xy}
*=<2cm,2cm>\txt{$u$}*\frm{-};
(14,0) *{[d]} ;
(-14,0) *{[a]} ;
(0,12) *{[c]} ;
(0,-12) *{[b]} ;
\end{xy}
\end{minipage}\]
for $[a], [b], [c], [d]\in G_1/H_1$.
The horizontal and vertical compositions are induced by those in $G$.

   The concept of quotient double groupoid using congruences was also studied in terms of quotient groupoids in \cite{11} (see also \cite{SD}).

\section*{Acknowledgement} We would like to thank to the referee for useful comments to improve the paper.

\end{document}